\documentstyle[12pt,amssymb,amsmath,amsfonts]{article}

\def \C{{\Bbb C}}

\def \Q{{\Bbb Q}}

\def \Z{{\Bbb Z}}
\def \ba{\mbox{{\boldmath $\alpha$}}}
\def \GQ{{\rm Gal}\,(\overline{{\Bbb Q}}/{\Bbb Q})}
\def \CP{{\mathcal P}}
\def \qed{\,\,\Box}

\def \phrase#1#2{

                   \noindent{\sc#1} {\sl#2}}
\def\begitem#1 {\bigskip\pagebreak[1]%
            \refstepcounter{subsection}{\nopagebreak[4]%
            \thesubsection\hskip 0.5truecm}
            {\sc#1}\hskip 1pt.\nopagebreak[4]\par\nopagebreak[4]%
             \begin{enumerate}\rm\nopagebreak[4]}
\def\BEGITEMK#1 #2{\bigskip\pagebreak[1]%
             \refstepcounter{subsection}{\nopagebreak[4]%
            \thesubsection\hskip 0.5truecm}\nopagebreak[4]
            {\bf#1}\hskip 1pt.\nopagebreak[4]\par\nopagebreak[4]%
            \medskip\nopagebreak[4]\rm#2\nopagebreak[4]%
            \begin{enumerate}\nopagebreak[4]\rm}
\def\enditem{\end{enumerate}}

\begin{document}
\title{Multivariate Diophantine equations\\with many solutions}
\author{J.-H. Evertse, P. Moree, C.L. Stewart
\footnote{The research of the
third author was
supported in part by Grant A3528 from the Natural Sciences and
Engineering Research
Council of Canada.}~ and R. Tijdeman}
\date{}
\maketitle
\vspace{0.5cm}

\begin{abstract}
\noindent
Among other things we show that for each 
$n$-tuple of positive rational numbers $(a_1, \ldots , a_n)$ there are 
sets of primes $S$ of arbitrarily large cardinality $s$
such that the solutions of the equation $a_1x_1+\cdots +a_nx_n=1$
with $x_1,\ldots,x_n$ $S$-units are not
contained in
fewer than $\exp ((4+o(1)) s^{1/2} (\log s)^{-1/2})$ proper linear subspaces 
of $\C^n$.
This generalizes a result of Erd\H{o}s, Stewart and Tijdeman
\cite{EST}  for $S$-unit equations in two variables.\\
Further, we prove that for any
algebraic number field $K$ of degree $n$, any integer $m$
with $1\leq m<n$, and any sufficiently large $s$
there are integers $\alpha_0,\ldots ,\alpha_m$
in $K$ which are linearly independent over $\Q$, and prime numbers
$p_1,\ldots ,p_s$, such that the norm polynomial equation
$|N_{K/\Q}(\alpha_0+\alpha_1x_1+\cdots +\alpha_mx_m)|=p_1^{z_1}\cdots
p_s^{z_s}$ has at least
$\exp \{ (1+o(1))\frac{n}{m}s^{m/n}(\log s)^{-1+m/n}\}$
solutions in $x_1,\ldots ,x_m,\,z_1,\ldots ,z_s\in\Z$.
This generalizes a result of Moree and Stewart \cite{MS}
for $m=1$.\\
Our main tool, also established in this paper, is an effective lower bound
for the number $\psi_{K,T}(X,Y)$
of ideals in a number field $K$ of norm $\leq X$
composed of prime ideals which lie outside a given
finite set of prime ideals $T$ and which have norm $\leq Y$. This generalizes
results of Canfield, Erd\H{o}s and Pomerance \cite{CEP} and of
Moree and Stewart \cite{MS}.
\\[0.5cm]
\end{abstract}
\newpage

\section{Introduction}
\setcounter{equation}{0}
\noindent
Let $S = \{ p_1,\ldots , p_s \}$ be a set of prime numbers. We call a 
rational
number an $S$-unit if both the denominator and the numerator of its
simplified representation are composed of primes from $S$.
Evertse \cite{Ev1} proved that for any non-zero rational numbers
$a,b$, the equation
$ax+by=1$ in $S$-units $x,y$ has at most
${\rm exp}(4s+6)$ solutions. On the other hand, Erd\H{o}s, Stewart and
Tijdeman \cite{EST} showed that equations of this type can have as many as
$\exp \{(4+o(1))(s/{\rm log}s)^{1/2} \}$ such solutions as
$s \to \infty$. Thus the dependence on $s$ cannot be polynomial.
In the present paper we generalise this result to $S$-unit equations in
an arbitrary number $n$ of variables. Here $n$ is considered to be given.
\\[0.5cm]
In \cite{Ev2} Evertse proved that for given non-zero rational numbers
$a_1,\ldots ,a_n$, the equation
\begin{equation}
a_1x_1+a_2x_2+\cdots +a_nx_n = 1\quad
\mbox{in $S$-units $x_1,x_2,\ldots ,x_n$}
\label{1}
\end{equation}
has at most $(2^{35}n^2)^{n^3(s+1)}$
non-degenerate solutions. We call a solution degenerate if there is some
non-empty proper subset $\{i_1,\ldots ,i_k \}$ of $\{ 1,\ldots ,n\}$ such 
that
$a_{i_1}x_{i_1}+a_{i_2}x_{i_2}+\cdots +a_{i_k}x_{i_k}=0$ and
otherwise non-degenerate.
In \cite{EGST}, Evertse, Gy\H{o}ry, Stewart and Tijdeman showed 
that there
are equations
(\ref{1}) which have as many as $\exp\{(4+o(1))(s/{\rm log}s)^{1/2} \}$
non-degenerate solutions as $s \to \infty$ and subsequently Granville 
\cite{Gr}
improved this to $\exp(c_0 s^{1-1/n} ( \log s)^{-1/n})$ for a positive 
number
$c_0$.
For our first result we shall establish a version of Granville's theorem 
with
$c_0$ given explicitly.
\\[0.3cm]
\phrase{Theorem 1. }
{Let $\varepsilon$ be a positive real number and let $a_1,\ldots
,a_n$ be non-zero
rational numbers.
There exists a positive number $s_0$, which is effectively computable in 
terms
of $\varepsilon$ and $a_1, \ldots , a_n$,
with the property that for every integer $s\geq s_0$
there is a set of primes
$S$ of cardinality $s$
such that
equation {\em (\ref{1})} has at least
\[
\exp
\big\{  (1-\varepsilon )\mbox{${\frac {n^2}{n-1}}$} 
s^{1-1/n}({\rm log}s)^{-1/n} \big\}
\]
non-degenerate solutions in $S$-units $x_1, x_2,\ldots ,x_n.$}
\\[0.5cm]
Theorem 1 does not exclude the possibility that the sets of solutions
of the equations
(\ref{1}) under consideration are of a special shape, for instance
that they are contained in the union
of a small number of proper linear subspaces of $\Q^n$ or in some
algebraic variety of small degree. We shall prove in Theorem 2 that
this is not the
case.
\\[0.3cm]
Let again $S$ be a set of primes and ${\bf a}=(a_1,\ldots ,a_n)$ a tuple
of non-zero rational numbers. 
Recall that the total degree of a polynomial $P$
is the maximum of the sums $k_1+\cdots +k_n$
taken over all monomials
$X_1^{k_1}\cdots X_n^{k_n}$ occurring in $P$.
Define $g({\bf a},S)$ to be the smallest 
integer
$g$ with the following property: there exists a polynomial
$P\in \C [X_1,\ldots ,X_n]$ of total degree
$g$, not divisible by
$a_1X_1+\cdots +a_nX_n-1$, such that
\begin{equation}
P(x_1,\ldots ,x_n)=0\quad\mbox{for every solution $(x_1,\ldots ,x_n)$
of (\ref{1}).}
\label{2}
\end{equation}
For instance, suppose that the set of solutions of (\ref{1})
is contained in the union of $t$ proper linear subspaces of $\C^n$,
given by equations $c_{i1}X_1+\cdots +c_{in}X_n=0$ $(i=1,\ldots ,t)$, say.
Then (\ref{2}) is satisfied by $P=\prod_{i=1}^t(\sum_{j=1}^n c_{ij}X_j)$
which is not divisible by $a_1X_1+\cdots +a_nX_n-1$;
hence $t\geq g({\bf a}, S)$. This means that if $g({\bf a}, S)$
is large, the set of solutions of (\ref{1}) cannot be contained in the
union of a small number of proper linear subspaces of $\C^n$. Likewise,
the set of solutions of (\ref{1}) cannot be contained in a proper algebraic
subvariety of small degree of the variety given by (\ref{1}).
Our precise result is as follows.
\\[0.3cm]
\phrase{Theorem 2. } {Let $\varepsilon$ be a positive real number and let
${\bf a} = (a_1, \ldots , a_n)$ be an
$n$-tuple of non-zero rational numbers.
There exists a positive number $s_1$, which is effectively computable
in terms of
$\varepsilon$ and ${\bf a}$, with the property that for every integer
$s\geq s_1$
there is a set of primes $S$ of cardinality $s$ such that}
\[
g({\bf a},S)\geq
\exp\big\{ (4-\varepsilon )s^{1/2}(\log s)^{-1/2}\big\}.
\]
\\[0.3cm]
Note that for $n=2$, both Theorems 1 and 2 imply the above mentioned result
of Erd\H{o}s, Stewart and Tijdeman.
\\[0.5cm]
We prove results analogous to Theorems 1 and 2 for ``norm polynomial
equations."
\\[0.2cm]
In what follows,
$K$ is an algebraic number field.
We denote by $O_K$ the ring of integers of $K$.
Let $\alpha_0,\ldots ,\alpha_m$ be elements of $O_K$ which are linearly
independent over $\Q$
and for which $\Q (\alpha_0,\ldots ,\alpha_m)=K$.
Further, let $p_1,\ldots ,p_s$ be distinct prime numbers.
From results of Schmidt \cite{WMS} and Schlickewei \cite{HPS}, it follows
that the
{\em norm form equation}
\begin{eqnarray}
|N_{K/\Q}(\alpha_0x_0+\cdots +\alpha_mx_m)|=p_1^{z_1}\cdots p_s^{z_s}
\label{30}\hspace*{5cm}\\
\mbox{in $x_0,\ldots ,x_m$, $z_1,\ldots ,z_s\in\Z$ with
$\gcd (x_0,\ldots ,x_m)=1$}\nonumber
\end{eqnarray}
has only finitely many solutions if and only if the
left-hand side satisfies some suitable non-degeneracy condition.
Instead of (\ref{30}) we deal with {\em norm polynomial equations}
\begin{equation}
|N_{K/\Q}(\alpha_0+\alpha_1x_1+\cdots +\alpha_mx_m)|=p_1^{z_1}\cdots 
p_s^{z_s}
\quad
\mbox{in $x_1,\ldots ,x_m$, $z_1,\ldots ,z_s\in\Z$,}\label{31}
\end{equation}
that is, norm form equations with $x_0=1$. As it turns out,
the number of solutions of a norm polynomial
equation is always finite, since there are no
degenerate cases. 
Recently B\'{e}rczes and Gy\H{o}ry (\cite[Theorem 2, Corollary 8]{BG}
or \cite[Chapter 1]{B})
proved a general result on
decomposable polynomial equations, a special case of which is that
equation (\ref{31}) has at most
\[
\big( 2^{17}n\big)^{\delta (m)(s+1)}
\]
solutions, where $n=[K: \Q ]$ and $\delta (m) =\frac{2}{3}(m+1)(m+2)(2m+3)-4$.
\\[0.3cm]
Note that for $m=1$,
equation (\ref{31}) is just the generalised Ramanujan-Nagell equation
\begin{equation}
|f(x)|=p_1^{z_1}\cdots p_s^{z_s}\quad\mbox{in $x,z_1,\ldots ,z_s\in\Z$,}
\label{32}
\end{equation}
where $f$ is an irreducible polynomial in $\Z [X]$ of degree at least $2$.
Erd\H{o}s, Stewart and Tijdeman \cite{EST} proved that given any $n\geq 2$
and any sufficiently large integer $s$
there are a polynomial $f\in\Z [X]$ of degree $n$
and primes $p_1,\ldots ,p_s$
such that (\ref{32})
has more than $\exp\{ (1+o(1))n^2s^{1/n}(\log s)^{(1/n)-1}\}$ solutions.
The polynomial constructed by
Erd\H{o}s,
Stewart and Tijdeman splits into linear factors over $\Q$.
\\[0.2cm]
Subsequently Moree and Stewart \cite{MS} proved a similar result
in which the constructed polynomial $f$ is irreducible.
More precisely,
let $K$  be a field of degree $n$ over $\mathbb{Q}$ and let $f$ be a monic
irreducible polynomial in $\mathbb{Z}$$[X]$ of degree $n$ such
that a root of
$f$ generates $K$ over $\mathbb{Q}$.
Let $\pi_f (x)$ denote the number of primes $p$ with $p \leq x$ for which
$f(x) \equiv 0\,$(mod $p$) has a solution.
It follows from the Chebotarev density theorem (see Theorems 1.3 and 1.4 of
\cite{LaOd}) that
\[
\pi_f (x) = \frac{1}{c_K} (1+o(1)) \frac{x}{\log x},
\]
where $c_K$ is a positive number which depends on $K$ only.
Let $L$ denote the normal closure of $K$. Then $c_K$ equals
$[L:\mathbb Q]$ divided by the
number of field automorphisms of $L/\mathbb Q$ that fix at least
one root of $f$, or in group theoretic
terms, $c_K=\# G/\#(\cup_{\sigma\in G}\sigma H\sigma^{-1})$,
where $H={\rm Gal}(L/K)$ and $G={\rm Gal}(L/\mathbb Q)$, see
\cite[Theorem 2]{BB}.
Thus $1 \leq c_K \leq n$ is a
rational number and if $K$ is normal then $c_K = n$.
Moree and Stewart \cite{MS} proved that for each field $K$ of degree $n$ 
over
$\mathbb{Q}$ there is a polynomial $f$, as above, such that the
number of solutions
of $(5)$ is $\exp \{ (1+ o(1)) n(c_Ks)^{1/n} (\log s)^{1/n -1} \}$.
\\[0.3cm]
\noindent
We generalize the result of Moree and Stewart to norm polynomial equations
as follows.
\\[0.3cm]
\phrase{Theorem 3. }
{Let $K$ be an algebraic number field of degree $n\geq 2$. Let
$\alpha_1,\ldots ,\alpha_m$ be elements of $O_K$ which are linearly
independent over $\Q$ where $1\leq m\leq n-1$. Let $\varepsilon >0$.
There exists a positive number $s_2$ which is effectively computable
in terms of
$\varepsilon$, $K$ and $\alpha_1, \ldots , \alpha_m$,
with the property that for every integer $s\geq s_2$
there are
a set
$S=\{ p_1,\ldots ,p_s\}$ of rational prime numbers and a number
$\alpha_0$ with
\begin{equation}
\alpha_0\in O_K,\,\, \Q (\alpha_0)=K,\,\,
\mbox{$\alpha_0$ $\Q$-linearly independent of $\alpha_1,\ldots ,\alpha_m$,}
\label{33}
\end{equation}
such that equation {\em (\ref{31})} has more than
\[
\exp\{ (1-\varepsilon )\mbox{$\frac{n}{m}$} (c_Ks)^{m/n}(\log s)^{(m/n)-1}\}
\]
solutions.}
\\[0.5cm]
Given a set of primes $S=\{ p_1,\ldots ,p_s\}$ and a tuple
$\ba =(\alpha_0,\ldots ,\alpha_n)$ of elements of $O_K$,
we define $g(\ba ,S)$ to be the smallest integer $g$ with
the following property: there exists a non-identically zero polynomial
$P\in\C [X_1,\ldots ,X_m]$ of total degree $g$ such that
\begin{equation}
P(x_1,\ldots ,x_m)=0\quad\mbox{for every solution
$(x_1,\ldots ,x_m,\, z_1,\ldots ,z_s)$ of (\ref{31}).}
\label{34}
\end{equation}
We prove the following result.
\\[0.3cm]
\phrase{Theorem 4. } {Let $K,n,m$, $\alpha_1,\ldots ,\alpha_m$ and
$\varepsilon >0$ be as in Theorem 3.
There exists a positive number $s_3$, which is effectively computable
in terms of
$\varepsilon, K$ and $\alpha_1 , \ldots , \alpha_m$,
with the property that for every integer $s\geq s_3$
there are
a set
$S=\{ p_1,\ldots ,p_s\}$ of rational prime numbers
and a number $\alpha_0$ with {\em (\ref{33})},
such that
\[
g(\ba ,S)
\geq\exp\{ (1-\varepsilon )n (c_Ks)^{1/n}(\log s)^{(1/n)-1}\}.
\]
Here $\ba =(\alpha_0,\alpha_1,\ldots ,\alpha_m)$.}
\\[0.5cm]
It should be noted that both Theorems 3 and 4 with $m=1$
imply the result of Moree and Stewart mentioned above.
\\[0.5cm]
The main tool in the proofs of Theorems 1-4 is a lower bound for the
number of ideals in a given number field which have norm $\leq X$,
are composed of prime ideals $\leq Y$, and which are composed of prime
ideals outside a given finite set of prime ideals $T$. 
We have stated this result below since it is not in the literature and since it 
may have some independent interest.
We first recall some history.
\\[0.2cm]
Let $\psi(X,Y)$ be the number of positive rational
integers not exceeding $X$
which are free of prime divisors larger than $Y$.
Canfield, Erd\H{o}s and Pomerance \cite{CEP} proved that there exists
an absolute constant $C$ such that if
$X,Y$ are positive reals with $Y\geq 3$
and with $u := \frac{\log X}{\log Y} \geq 3$, then
\begin{eqnarray}\label{CEP}
\psi(X,Y) \geq X\exp\left\{-u\Big\{\log (u\log u)-1+\frac{\log_2u-1}{\log u}
+C\cdot\Big(\frac {\log_2 u}{\log u}\Big)^2\Big\}\right\},
\end{eqnarray}
where $\log_2 u =\log\log u$.
Further, 
Hildebrand \cite{H86} showed that for arbitrary fixed $\varepsilon >0$,
one has uniformly under the condition $X\ge 2,~\exp\{(\log_2 X)^{{5\over
3}+\varepsilon}\}\le Y\le X$,
\begin{equation}
\label{hildesharp}
\psi(X,Y)=X\rho(u)\Big\{1+O\left({\log(u+1)\over \log Y}\right)\Big\},
\end{equation}
where $\rho(u)$
denotes the Dickman-de Bruijn function.
\\[0.2cm]
More generally, let $K$ be a number field.
By an ideal of the ring of integers $O_K$ we shall mean a non-zero ideal.
Denote by $\psi_K(X,Y)$ the number of ideals of $O_K$ with
norm at most $X$ composed of prime ideals of $O_K$ of norm at most $Y$.
Here the norm of an ideal ${\frak a}$ is the cardinality of the residue
class ring $O_K/{\frak a}$.
By Moree and Stewart \cite[Theorem 2]{MS} 
there exists a constant $C_K>0$,
depending only on $K$,
such that with $X,Y$ and $u$ as above we have
\begin{eqnarray}\label{MS}
\psi_K(X,Y) \geq X\exp\left\{-u\Big\{\log (u\log u) -1
+\frac{\log_2u-1}{\log u}+C_K\Big(\frac{\log_2u}{\log u}\Big)^2\Big\}
\right\}.
\end{eqnarray}
This result has been proved by extending the method of Canfield,
Erd\H{o}s and Pomerance.
\\[0.2cm]
Now let $T$ be a finite set of prime ideals of $O_K$, 
and denote by $\psi_{K,T}(X,Y)$
the number of ideals of $O_K$ which have norm $\leq X$ and are composed
of prime ideals which have norm $\leq Y$ and lie outside $T$.
We prove the following:
\\[0.3cm]
\phrase{Theorem 5. } {There exists a positive
effectively computable number $C_{K,T}$ depending only on $K$ and $T$ such
that for $X,Y\ge 1$ with $u:=\frac{\log X}{\log Y}\ge 3$ we have
\begin{eqnarray}
\psi_{K,T}(X,Y) &\geq& X\exp\Big\{-u\big\{\log (u\log u)-1\, +
\nonumber\\[0.1cm]
&&
\label{psi}
\qquad\qquad\qquad\qquad
+\,\frac{\log_2 u -1}{\log u}+C_{K,T}\Big(\frac{\log_2 u}{\log u}\Big)^2\big\}
\Big\}.
\end{eqnarray} 
}
\\[0.2cm]
In the proof of Theorem 5 we did not use the ideas of
Canfield, Erd\H{o}s and Pomerance, but instead extended the arguments
from Hildebrand's paper \cite{H86} mentioned above. 
Another more straightforward method
to obtain a lower bound for $\psi_{K,T}$ such as (\ref{psi})
is by combining the estimate
(\ref{MS}) for $\psi_K(X,Y)$ with an interval result for $\psi_K(X,Y)$
due to Moree \cite{M3}. Unfortunately, the result obtained by this approach
is valid only for a much smaller $X,Y$-range,  
and it is not at all transparent whether the constant
$C_{K,T}$ arising from this approach is effective. 
In \cite{BH} Buchmann and
Hollinger, assuming the
Generalized Riemann Hypothesis, established
a non-trivial lower bound for $\psi_K(X,Y)$, uniform in $K$,
involving the
degree of the normal closure and the discriminant $D_K$ of $K$. They
did so by
using
the method of Canfield, Erd\H{o}s and Pomerance. Our method
to prove Theorem 5 can be used to obtain a variant of the result of
Buchmann and Hollinger
with much smaller error term. 
As a starting point in our approach one may
take equation (11.RH) of Lang \cite{LAnalytic}.

\section{Proof of Theorem 5.}
\setcounter{equation}{0}
\noindent 
We recall some properties of the Dickman-de
Bruijn function $\rho (u)$. This function is the unique
continuous solution of the differential-difference equation
$u\rho'(u)=-\rho(u-1)$ for $u>1$ with initial condition
$\rho(u)=1$ in the interval $[0,1]$
(and, by
convention, $\rho(u):=0$ for $u<0$).
Recall that according to Hildebrand's estimate (\ref{hildesharp}),  
$\rho(u)$ is the density of the set of integers $\leq X$ composed of
prime numbers $\leq X^{1/u}$ as $X$ tends to infinity;
therefore, $0\le \rho(u)\le 1$.
In the following lemma we have collected some further easily provable
properties of the Dickman-de Bruijn function that
will be needed in the sequel.
\\[0.3cm]
\phrase{Lemma 1. } {We have\\
i) $u\rho(u)=\int_{u-1}^u \rho(t)dt~~{\rm for~}u\ge 1.$\\
ii) $\rho(u)>0$ for $u>0$.\\
iii) $\rho(u)$ is decreasing for $u>1$.\\
iv) $-\rho'(u)/\rho(u)$ is increasing for $u>1$.\\
v) $-\rho'(u)\le \rho(u)\log(2u\log^2(u+3))$ for $u>0$, $u\ne 1$.\\
vi) $\rho(u-t)\le \rho(u)4(2u\log^2(u+3))^t$
for $u\ge 0$ and $0\le t\le 1$.}
\\[0.3cm]
{\bf Proof}. This is in essence
\cite[Lemma 6]{H86}, see also
\cite[p. 30]{MProef}. Parts v) and vi) are, however, modified so as
to obtain explicit estimates valid for $u>0$.
They require some easy numerical verifications that
are left to the interested reader.$\qed$
\\[0.5cm]
An important quantity in the study of the Dickman-de Bruijn function
is the function $\xi(u)$. For any given $u>1$, $\xi(u)$ is defined
as the unique positive solution of the transcendental equation
\begin{equation}
\label{transcendental}
{e^{\xi}-1\over \xi}=u.
\end{equation}
The quantity $\xi (u)$ exists and is unique, since 
$\lim_{x\downarrow 0} (e^x-1)/x =1$ and since
$(e^x-1)/x$ 
is strictly increasing for $x>0$.
The Fourier 
transform ${\hat \rho}$ of $\rho$ involves $\xi$. By
writing $\rho$ as the Fourier transform of ${\hat \rho}$ and
applying the saddle point method one obtains \cite{A} that for $u\ge 1$,
\begin{equation}
\label{alladi}
\rho(u)=\sqrt{\xi'(u)\over 2\pi}\exp\Big\{\gamma-\int_1^u\xi(t)dt\Big\}
\{1+O({1\over u})\}.
\end{equation}
(It is not difficult to show that $\xi'(u)\sim 1/u$ as $u$ tends to
infinity.)
For our purposes we need an effective lower bound of
the quality of (\ref{alladi}). The next lemma fulfils our needs.
\\[0.3cm]
\phrase{Lemma 2. } {For $u\ge 1$ we have
$$\exp\Big\{-\int_2^{u+1}\xi(t)dt\Big\}\le \rho(u)\le 
\exp\Big\{-\int_1^{u}\xi(t)dt\Big\}.$$
}
\\[0.2cm]
{\bf Proof}. Let
$f(u)=-\rho'(u)/\rho(u)$
denote the logarithmic derivative of $1/\rho(u)$.
Using parts i) and iv) of Lemma 1 we deduce that
$$u=\int_{u-1}^u{\rho(t)\over \rho(u)}dt=\int_{u-1}^u
e^{\int_t^uf(s)ds}dt\le \int_{u-1}^ue^{(u-t)f(u)}dt
={e^{f(u)}-1\over f(u)},$$
and thus, by the monotonicity of $(e^x-1)/x$, that $f(u)\ge \xi(u)$
for $u>1$. By a similar argument we find that $f(u)\le \xi(u+1)$
for $u>0$ and $u\ne 1$.
On noting that
\[
\exp\Big(-\int_{1}^u\xi(s+1)ds\Big)\le 
\rho(u)=\exp\Big(-\int_{1}^uf(s)ds\Big)\le
\exp\Big(-\int_{1}^u\xi(s)ds\Big),
\]
the proof is completed.$\qed$
\\[0.5cm]
The method of bootstrapping allows one to obtain an asymptotic expression for
$\xi(u)$ with error $O(\log^{-k}u)$ for arbitrarily large $k$.
To illustrate this we do the first few iterations.
From (\ref{transcendental}) we deduce that
\begin{equation}
\label{eersteronde}
\xi=\log \xi + \log u + O({1\over \xi\cdot u}),~\xi\cdot u\rightarrow 
\infty.
\end{equation}
Notice that for $u$ sufficiently large $1<\xi<2\log u$. It follows from
(\ref{eersteronde}) that $\xi=\log u+ O(\log_2 u)$. Substituting this into
the right-hand side of
(\ref{eersteronde}) then yields
$\xi=\log u+\log_2 u+O({\log_2 u/\log u})$. Note that the implied
constant is effective. By repeatedly
substituting
the lastly found asymptotic expression for $\xi (u)$ into the right-hand side
of (\ref{eersteronde}), one can calculate an asymptotic expression for
$\xi (u)$ with error $O(\log^{-k}u)$ for arbitrary $k>1$, with
effective implied constant. This then implies, using
Lemma 2, that for arbitrary $k>1$ we can find an elementary
explicit function $g_k(u)$ such that
$\rho(u)\ge \exp(g_k(u)+O_k(u\log^{-k}u))$, where the implied constant
is effective.
For example, by substituting $\xi=\log u+\log_2 u+O({\log_2 u/\log u})$
into the right-hand side of (\ref{eersteronde}) we obtain for $u\geq 3$,    
\[
\xi =\log u +\log_2 u +\frac{\log_2 u}{\log u}+
O\Big(\big(\frac{\log_2u}{\log u}\big)^2\Big).
\]
Using Lemma 2 we then find that, for $u\geq 3$,
\begin{equation}
\label{rhoexpansie}
\rho(u)\ge \exp\left\{-u\Big\{\log (u\log u)-1+{\log_2 u-1\over \log u}
+O\Big(\big({\log_2 u\over \log u}\big)^2\Big)\Big\}\right\},
\end{equation}
where the implied constant is effective.
\\
Alternatively $g_k(u)$ can be computed using the convergent series expansion
for $\xi(u)$ on pp. 145, 146 of Hildebrand and Tenenbaum \cite{HT}.
\\[0.2cm]
Now let $K$ be an algebraic number field. 
We put $P({\frak a})={\rm max}\{N{\frak p}~:~{\frak p}|{\frak a}\}$ 
for an ideal ${\frak a}\ne (1)$ of $O_K$
and
$P((1))=1$ (here and in the sequel the symbol $\frak p$ is exclusively
used to indicate a prime ideal).
We denote by $N_K(Y)$ the number of ideals of $O_K$ of norm $\leq Y$ and
for a given finite set of prime ideals $T$ of $O_K$, by $N_{K,T}(Y)$
the number of ideals of $O_K$ of norm $\leq Y$ which are coprime with
each of the prime ideals from $T$. For instance from
the arguments in Lang \cite[Chap. VI-VIII]{LBoek} 
it follows that
\[
N_K(Y)=A_KY+O(Y^{1-1/[K:\Bbb Q]}),
\quad\mbox{where $A_K={\rm Res}_{s=1}\zeta_K(s)$}
\]
is the residue of the Dedekind
zeta-function at $s=1$ (which as is well-known can be expressed in terms
of invariants such as the class number and regulator of the field $K$)
and where the implied constant is effective and depends only on $K$.
By means of the principle of inclusion and exclusion it then follows that
\begin{equation}
\label{ideaalafschatting}
N_{K,T}(Y)=A_{K,T}Y+O(Y^{1-1/[K:\Bbb Q]}),
\quad\mbox{where $A_{K,T}=A_K\prod_{{\frak p}\in T}(1-{1\over N{\frak p}})$}
\end{equation}
and where the implied constant is effective and depends only on $K$ and $T$.
\\[0.2cm]
As before, we 
denote by $\psi_{K,T}(X,Y)$
the number of ideals of $O_K$ of norm at most $X$ which are composed
of prime ideals which do not belong to the finite set of prime ideals $T$ 
and, moreover, have norm at
most $Y$. 
The ideals so counted form a free arithmetical semigroup
satisfying the conditions of Theorem 1 of  \cite[Chapter 4]{MProef}.
It then follows that, for
arbitrary fixed $\varepsilon\in (0,1)$, we have uniformly for
$1\le u\le (1-\varepsilon)\log_2 X/\log_3 X$ that
\begin{equation}
\label{asymptotiek}
\psi_{K,T}(X,Y)\sim A_{K,T}X\rho(u)\quad\mbox{as $X\to\infty$,}
\end{equation}
where $\log_3X=\log\log\log X$.
Thus we get a density interpretation of $\rho(u)$ similar
to that for $\psi(X,Y)$.\\
\\[0.2cm]
The proof of (\ref{asymptotiek}) is based on the Buchstab 
functional equation for free
arithmetical semigroups. In order to obtain Theorem
5, which gives a lower bound
for $\psi_{K,T}(X,Y)$ valid for a much larger $X,Y$-region, a different
functional equation will
be used. This equation along
with several other ideas that go into the proof of Theorem 5 are due
to Hildebrand \cite{H86}, cf. \cite[pp. 388-389]{T}, who worked in the case 
where $K=\mathbb Q$
and $T$ is the empty set.
Put ${\frak q}=\prod_{{\frak p}\in T}{\frak p}$.
Define
\[
\Lambda_{K,T}(\frak a)=
\left\{
\begin{array}{l}
\log N\frak p\quad\mbox{if $\frak a=\frak p^m$ for some
${\frak p}\not\in T$ and $m\geq 1$,}\\[0.2cm]
0\quad\mbox{otherwise.}
\end{array}\right.
\]
Then for $X\geq Y$ we have
\begin{eqnarray}\label{func}
\psi_{K,T}(X,Y)\log X=\int_1^X{\psi_{K,T}(t,Y)\over t}dt
+\sum_{N{\frak a}\le X\atop P({\frak a})\le Y}\Lambda_{K,T}({\frak a})
\cdot \mbox{$\psi_{K,T}({X\over N{\frak a}},Y)$}.
\end{eqnarray}
In order to establish the validity
of this equation we express the sum of all $\log N\frak a$
with $\frak a$ satisfying $N{\frak a}\le X$,  $P({\frak a})\le Y$
and ${\frak a}$ coprime with ${\frak q}$ in two different ways.
On the one hand we find by integration by parts that this sum can
be expressed as
\[
\psi_{K,T}(X,Y)\log X-\int_1^X{\psi_{K,T}(t,Y)\over t}dt,
\]
on the other hand we notice that the sum can be rewritten as follows
\[
\begin{array}{rl}
\displaystyle{
\sum_{{N{\frak a}\le X,~{\frak a}+{\frak q}=(1)}\atop P({\frak a})\le Y}
\sum_{{\frak b}|{\frak a}}\Lambda_{K,T}({\frak b})}&
=\displaystyle{
\sum_{N{\frak b}\le X\atop P({\frak b})\le Y}\Lambda_{K,T}({\frak b}) 
\sum_{{N{\frak a}\le X,~
{\frak a}+{\frak q}=(1)}\atop {\frak b}|{\frak a},~P({\frak a})\le Y}1}
\\ 
&\quad
\\
&=\displaystyle{\sum_{N{\frak b}\le X\atop P({\frak b})\le Y}
\Lambda_{K,T}({\frak b})
~\psi_{K,T}( \mbox{${X\over N{\frak b}}$},Y),}
\end{array}
\]
where we used that 
$\log N{\frak a}=\sum_{{\frak b}|{\frak a}}\Lambda_{K,T} ({\frak b})$ 
for any
ideal $\frak a$ coprime with ${\frak q}$. 
\\
Using functional equation (\ref{func})
and Lemmata 3 and 4 below, we will deduce the crucial Lemma 5,
and from that, Theorem 5. 
\\[0.3cm]
\phrase{Lemma 3. } {Let $K$ be a number field and $T$ a finite set
of prime ideals in $O_K$. Put
$\log^+ Y=\max\{1,\log Y\}$. Then
\[
\sum_{N{\frak a}\le Y}{\Lambda_{K,T}({\frak a})\over N{\frak a}}
=\log Y+c_{1,K,T}+E(Y)\quad
\mbox{for $Y\geq 1$,}
\]
where $c_{1,K,T}$ is a constant depending on $K$ and $T$ and where
for every $m\geq 1$ we have $|E(Y)|\le c'_m(\log ^+ Y)^{-m}$,
with $c'_m$ an effectively computable constant depending
on $m$, $K$ and $T$.}
\\[0.2cm]
{\bf Proof}. Let $\Pi_K(Y)$ denote the number of prime ideals of $K$
of norm $\leq Y$. Theorems 1.3, 1.4 of
Lagarias and Odlyzko \cite{LaOd}
imply an effective version
of the Prime Ideal Theorem of the shape $\Pi_K(Y)=Li(Y)+E_0(Y)$
where $Li(Y)=\int_2^Y (\log t)^{-1}dt$ and 
$|E_0(Y)|\leq c''_mY(\log^+ Y)^{-m}$ for every $m\geq 2$, with
$c''_m$ an effectively computable constant depending on $m$ and $K$. 
From this and the standard Stieltjes integration
and partial summation arguments one obtains Lemma 3.$\qed$
\\[0.3cm]
\phrase{Lemma 4.} {Let $0<\theta\le 1$, $m\ge 4$,  $1\le u\le Y^2$,
$Y\ge e^{m^{3m}}$ and let
$c'_m$ be as in Lemma 3.
Put
\[
S_{\theta}=\sum_{N{\frak a}\le Y^{\theta}}{\Lambda_{K,T}({\frak a})\over 
N{\frak a}}\cdot\mbox{$\rho (u-{\log N{\frak a}\over \log Y})$.}
\]
Then
\[
S_{\theta}=\log Y\int_0^{\theta}\rho(u-v)dv+E_1(\theta),
\]
with
\[
|E_1(\theta)|\le 17c'_m\rho(u)\Big\{2+{u\log^2(u+3)\over \log ^{m-1} 
Y}\theta^{-m}\Big\}.
\]
}
\\[0.2cm]
{\bf Proof}.
Using Lemma 3 we find by Stieltjes integration that
$$S_{\theta}=\int_0^{\theta}\rho(u-v)
d\left(\sum_{N{\frak a}\le Y^v}{\Lambda_{K,T}({\frak a})\over N
{\frak a}}\right)
=
\log Y\int_0^{\theta}\rho(u-v)dv+I_1(\theta)+I_2(\theta),$$
where
$I_1(\theta)=E(Y^{\theta})\rho(u-\theta)-E(1)\rho(u)$ and
$I_2(\theta)=\int_0^{\theta}\rho'(u-v)E(Y^v)dv$. Using Lemma 1 vi)
we
deduce that
$$|I_1(\theta)|\le c'_m\rho(u)\Big\{1+{8u\log^2 (u+3)\over \log^m 
Y}\theta^{-m}\Big\}.$$
For notational convenience let us put 
$g(u):=\log(2u\log^2 (u+3))$.
Then using Lemma 1 v), vi) we obtain
$$|I_2(\theta)|\le 4\rho(u)g(u)
\Big\{c'_m\int_0^{\log^{-1} Y}e^{vg(u)}dv+\int_{\log^{-1} 
Y}^{\theta}e^{vg(u)}|E(Y^v)|dv\Big\}.$$
The conditions on $u$
and $Y$ ensure that the first integral in the latter estimate is bounded 
above
by $g(u)^{-1}\exp(g(u)/\log Y)\le 8/g(u)$.
We split up the integration range of
the second integral at $\theta\log^{-1/m}Y$ and denote the corresponding 
integrals
by $I_3(\theta)$ and $I_4(\theta)$, respectively.
We have
\begin{equation}\label{i3}
|I_3(\theta)|\le c'_m{e^{\theta g(u)\log^{-{1\over m}}Y}\over \log^m Y}
\int_{\log^{-1}Y}^{\theta\log^{-{1\over m}}Y}{dv\over v^m}\le
{c'_m\over \log Y}e^{\theta g(u)/\log^{1\over m}Y}
\end{equation}
and
\begin{equation}\label{i4}
|I_4(\theta)|\le {c'_m\theta^{-m}\over 
\log^{m-1}Y}\int_{\theta\log^{-{1\over m}}Y}
^{\theta}e^{vg(u)}dv\le {c'_m\theta^{-m}\over \log^{m-1}Y}
{2u\log^2(u+3)\over g(u)}.
\end{equation}
Note that if $g(u)\le \log^{1/m}Y$, then $g(u)|I_3(\theta)|\le c'_m/4$.
If $g(u)>\log^{1/m}Y$, then 
thanks to our assumption $Y\ge e^{m^{3m}}$,
the right-hand side of (\ref{i3}) is smaller than the right-hand side of
(\ref{i4}), therefore both $|I_3(\theta )|$ and $|I_4(\theta )|$ are bounded
above by 
${c'_m\theta^{-m}\over \log^{m-1}Y}
{2u\log^2(u+3)\over g(u)}$. 
On adding the various estimates, our lemma follows.$\qed$
\\[0.3cm]
\phrase{Lemma 5. } {Let
$m\ge 4$ be arbitrary and
$1\le u\le Y^2$. Suppose that $Y\ge \max\{e^{m^{3m}},e^{1500c'_m}\}$.
Then
$$\psi_{K,T}(X,Y)\ge X\rho(u)\Delta\exp\Big(-1224c'_m\Big\{{\log(6(u+1))\over \log Y}+
{5\cdot 2^{m-1}(u+1)\over \log^{m-3} Y}\Big\}\Big),$$
where $\Delta:=\inf_{Y\ge 1}N_{K,T}(Y)/Y$.}
\\[0.3cm]
{\bf Proof}.
We set
$\delta(u):=\inf_{0\le v\le u} \psi_{K,T}(Y^v,Y)/(Y^v\rho(v))$.
Note that
$\delta(u)\ge \Delta$ for $0\le u\le 1$.
Let $u>1$.
Functional equation (\ref{func}) gives rise to the estimate
\begin{eqnarray*}
\psi_{K,T}(X,Y)\log X &\geq& \sum_{N{\frak a}\le Y}
\mbox{$\Lambda_{K,T}({\frak a})\psi_{K,T}({X\over N{\frak a}},Y)$}\\[0.2cm]
&\geq&
\mbox{$X\delta(u)S_{1\over 2}+X\delta(u-{1\over 2})(S_1-S_{1\over 2})$}.
\end{eqnarray*}
By dividing this inequality by $X\rho(u)\log X=Xu\rho(u)\log Y$
and then using Lemma 4, Lemma 1 i) and the fact that $\delta$ is decreasing,
we obtain
\[
{\psi_{K,T}(X,Y)\over X\rho(u)} \geq 
\mbox{$\delta(u)r(u)+
\delta(u-{1\over 2})\{1-r(u)-2|E_1(\mbox{$\frac{1}{2}$})|-|E_1(1)|\}$,}
\]
where
\[
r(u)={1\over u\rho(u)}\int_0^{{1\over 2}}\rho(u-v)dv\, .
\]
Since by Lemma 1 iii), $\rho$ is decreasing it follows that 
$r(u)\leq \frac{1}{2}$. Further,
\[
2|E_1(\mbox{$\frac{1}{2}$})|+|E_1(1)|
\leq f_m(u):= {51c'_m\over \log Y}\Big\{{2\over u}+
{5\cdot 2^m\over \log^{m-3} Y}
\Big\}.
\]
Hence 
\begin{eqnarray}
\label{deltarho}
{\psi_{K,T}(X,Y)\over X\rho(u)} \geq 
\mbox{$\frac{1}{2}\delta (u) +(\frac{1}{2}-f_m(u))\delta (u-\frac{1}{2})$.}
\end{eqnarray}
We want to establish that
\begin{equation}
\label{shiftback}
\mbox{$\delta(u)\ge 
\min (\Delta ,\delta(u-{1\over 2}))e^{-6f_m(u-{1\over 2})}$.}
\end{equation}
If $\delta(u)=\delta(u-{1\over 2})$, this inequality is trivially true.
If $\delta (u)=\delta (1)$ the inequality is true as well,
since $\delta (1)\geq \Delta$.
So assume that $\delta(u)<\delta(u-{1\over 2})$ and $\delta (u)<\delta (1)$.
Choose $\varepsilon$ with $0<\varepsilon <1$.
Then there exists $u'\in (\max (1,u-{1\over 2}),u]$ such that
$\psi_{K,T}(X',Y)/(X'\rho (u'))\le \delta(u)(1+\varepsilon)$, with
$X'=Y^{u'}$.
Using (\ref{deltarho}) with $u'$ replacing $u$ we then infer
\[
\begin{array}{rl}
\delta (u)(1+\varepsilon)&
\geq \frac{1}{2}\delta (u')+ (\frac{1}{2}-f_m(u'))\delta (u'-\frac{1}{2})
\\[0.1cm]
&\geq \frac{1}{2}\delta (u)+ (\frac{1}{2}-f_m(u-\frac{1}{2}))\delta (u-\frac{1}{2}).
\end{array}
\]
Since $\varepsilon$ may be chosen arbitrarily small,
the latter inequality implies that
$\delta(u)\ge \delta(u-{1\over 2})(1-2f_m(u-{1\over 2}))$.
The lower bound $Y\ge \exp(1500c'_m)$ ensures that $f_m(u-{1\over 2})<1/6$
and hence the validity of (\ref{shiftback}).
\\[0.2cm]
We now
iterate (\ref{shiftback}), the last step being with an argument $u_0>1$
such that $\delta (u_0-\frac{1}{2})\geq\Delta$. 
Since $f_m$ is decreasing,
this yields
$\delta(u)\geq \Delta \exp\{-6\sum_{k=0}^{2[u]}f_m({k+1\over 2})\}$.
Then Lemma 5 follows after an easy computation.$\qed$
\\[0.5cm]
{\bf Proof of Theorem 5.}
By (\ref{ideaalafschatting}) (which is effective),
there is
an effective constant $\Delta_0$ such that $\Delta\geq \Delta_0>0$.
Now from this fact, Lemma 5 with $m=6$ and (\ref{rhoexpansie}) (where the
implied constant can be made effective)
we obtain (\ref{psi})
with some effective constant $C_{K,T}>0$,
provided that $1\le u\le Y^2$ and $Y\geq Y_0$, where $Y_0$ is some effectively
computable number depending on $K$ and $T$.  
Note that if $u>Y^2$ and $Y\geq Y_1$ 
(with $Y_1\geq Y_0$ effective and depending
on $K,T$ and $C_{K,T}$)
then the right-hand side of (\ref{psi}) is
$<1$
so that (\ref{psi}) is trivially true
(as $\psi_{K,T}(X,Y)\geq 1$). 
Further, if
$Y\leq Y_1$ then for $X$ exceeding some effectively computable number $X_0$
depending on $K,T, Y_1$ and $C_{K,T}$ we have again that the
right-hand side of (\ref{psi}) is $<1$, so that (\ref{psi}) holds.
We can achieve that (\ref{psi}) holds
for the remaining values for $X,Y$, i.e., $Y\leq Y_1$ and $X\leq X_0$,
by enlarging the constant $C_{K,T}$ if necessary.
This completes the proof of Theorem 5.$\qed$
\\[0.3cm]
{\bf Remark.} Given any
$k>0$, a refinement of Theorem 5 with
error term
$\exp\{O(u\log^{-k}u)\}$ and effective implied constant can be given
by carrying out the bootstrap process for $\xi(u)$ far enough.

\section{Preparations for the proofs of Theorems 1--4.}
\setcounter{equation}{0}
\noindent
We start with a simple result on polynomial equations.
\\[0.3cm]
\phrase{Lemma 6. } {Let $Q\in \C [X_1,\ldots ,X_m]$ be a
non-trivial polynomial of total
degree $g$. Let $A,B\in \Z$ with $A<B$. Then the set of vectors
${\bf x}=(x_1,\ldots ,x_m)\in\Z^m$ with
\begin{equation}
Q({\bf x})=0,\qquad A\leq x_i\leq B\,\,\mbox{for $i=1,\ldots ,m$}
\label{3}
\end{equation}
has cardinality at most $g\cdot (B-A+1)^{m-1}$.}
\\[0.3cm]
{\bf Proof.} We proceed by induction on $m$. For $m=1$ the lemma is obvious.
Suppose $m>1$. Assume the lemma holds true for polynomials in fewer than $m$
variables. We may write
\[
Q(X_1,\ldots , X_m)=\sum_{i=0}^h Q_i(X_1,\ldots ,X_{m-1})X_m^i
\]
with $h\leq g$, $\, Q_i\in \C [X_1,\ldots ,X_{m-1}]$ of total degree $\leq 
g-i$
for $i=0,\ldots ,h$
and with $Q_h$ not identically zero.
Let $V$ be the set of tuples ${\bf x}$ with (\ref{3}).
Given ${\bf x}=(x_1,\ldots ,x_m)\in V$ we write
${\bf x'}=(x_1,\ldots ,x_{m-1})$.
\\[0.3cm]
First consider those ${\bf x}\in V$ for which $Q_h({\bf x'})\not=0$.
There are at most $(B-A+1)^{m-1}$ possibilities for ${\bf x'}$.
Fix one of those ${\bf x'}$. Substituting $x_i$ for $X_i$ $(i=1,\ldots 
,m-1)$
in $Q$ gives a non-zero polynomial of degree $h$ in $X_m$. Hence for given
${\bf x'}$ there are at most $h$ possibilities for $x_m$ such that
$Q({\bf x})=0$. So altogether, there are at most $h(B-A+1)^{m-1}$
vectors ${\bf x}\in V$ with $Q_h({\bf x'})\not=0$.
\\[0.3cm]
Now consider those ${\bf x}\in V$ for which $Q_h({\bf x'})=0$.
Recall that $Q_h$ has total degree at most $g-h$. So by the induction
hypothesis, there are at most $(g-h)(B-A+1)^{m-2}$ possibilities for
${\bf x'}$. For a fixed ${\bf x'}$, there are at most $B-A+1$
possibilities for $x_m$. Therefore,
there are at most $(g-h)(B-A+1)^{m-1}$
vectors ${\bf x}\in V$ with $Q_h({\bf x'})=0$.
\\[0.3cm]
Combining this with the upper bound $h(B-A+1)^{m-1}$ for the number of
vectors in $V$ with $Q_h({\bf x'})\not=0$,
we obtain that $V$ has cardinality at most
$g(B-A+1)^{m-1}$. This proves Lemma 6.$\qed$
\\[0.5cm]
Let $K$ be a number field.
We denote by $\xi\mapsto \xi^{(i)}$ $(i=1,\ldots ,[K:\Q ])$ the
isomorphic embeddings of $K$ into $\C$.
The prime ideal decomposition
of $\alpha\in O_K$ is by definition the prime ideal decomposition
of the principal ideal $(\alpha )$ generated by $\alpha$. We say that
$\alpha\in O_K$ is coprime with the ideal ${\frak a}$ if
$(\alpha )+{\frak a}=(1)$. 
\\[0.3cm]
\phrase{Lemma 7. } {Let $[K:\Q ]=n$.
Let ${\frak a}$ be an ideal of $O_K$ and let
$\alpha\in O_K$ be coprime to ${\frak a}$. 
Further, let $T$ be the set of prime ideals dividing ${\frak a}$.
Then there are effectively computable constants $C_1,C_2,C_3>1$, 
depending only on
$K,{\frak a}$ such that
for $X,Y$ with
$X>Y\geq C_1$,
the number of non-zero $\xi\in O_K$ with
\begin{equation}
\left.
\begin{array}{l}
\mbox{$|\xi^{(i)}|\leq C_2X^{1/n}$ for $i=1,\ldots ,n$},\\
\xi\equiv\alpha\pmod{{\frak a}},\\
\mbox{$(\xi)$ is composed of prime ideals of norm $\leq Y$\quad\quad}
\end{array}
\right\}
\label{22}
\end{equation}
is at least $C_3^{-1}\psi_{K,T}(X,Y)$.}
\\[0.3cm]
{\bf Proof.}
Below,
constants implied by $\ll$, $\gg$
depend only on $K$,${\frak a}$ and are all effective.
For $\xi\in O_K$ let $||\xi ||$
denote the maximum of the absolute values of the conjugates of $\xi$.
Denote by $h$ the class number of $K$.
By the effective version of the Chebotarev density theorem
from \cite{LaOd} (Theorems 1.3, 1.4)
each ideal class of $K$ contains a prime ideal outside $T$ 
with norm bounded above effectively in terms of $K,{\frak a}$.
Let $\mathcal{H}$ consist of one such prime ideal from each ideal class.
\\[0.2cm]
Assume that $Y$ exceeds the norms of the prime ideals from $\mathcal{H}$.
Let ${\frak b}$ be an ideal of norm at most $X$ composed of prime ideals
of norm at most $Y$ lying outside $T$.
Choose ${\frak p}$
from $\mathcal{H}$ such that ${\frak b}\cdot {\frak p}$ is a principal
ideal, $(\beta )$, say. Then $(\beta )$ has norm $\ll X$
and is composed of prime ideals of norm $\leq Y$ lying outside $T$.
Further,
there are at most $h$ ways of obtaining a given principal ideal $(\beta )$
by multiplying an ideal of norm at most $X$ 
with a prime ideal from $\mathcal{H}$.
Therefore,
the number of principal ideals of norm $\ll X$, composed
of prime ideals of norm at most $Y$ and lying outside $T$, is
at least $h^{-1}\psi_{K,T}(X,Y)$.
\\[0.2cm]
We choose from each residue class in $(O_K/{\frak a})^*$
a representative $\gamma$ for which
$||\gamma ||$ is minimal. Denote the set of these representatives by
$\mathcal{R}$. Suppose $\mathcal{R}$ has cardinality $m$. 
Clearly,
each element from $\mathcal{R}$ is composed of prime ideals outside $T$.
Furthermore, of each element of $\mathcal{R}$ 
the absolute value of the norm can be bounded
above effectively in terms of $K,{\frak a}$.
\\[0.2cm]
Assume that $Y$ exceeds the absolute values of the
norms of the elements from $\mathcal{R}$.
Then the elements of $\mathcal{R}$ are composed of prime ideals outside $T$
of norm at most $Y$.
Take a principal ideal $(\beta )$ of norm $\ll X$ composed of prime ideals
of norm at most $Y$ lying outside $T$.
According to, for instance, \cite{ST}, Lemma A.15, there is a $\beta '$
with $(\beta ')=(\beta )$ and $||\beta ' ||\ll X^{1/n}$.
Clearly, $\beta '$ is coprime with ${\frak a}$, so there is a
$\gamma\in\mathcal{R}$ with $\xi :=\beta '\gamma\equiv\alpha\pmod{{\frak a}}$.
Note that $||\xi ||\ll X^{1/n}$, and that $(\xi )$ is composed of prime
ideals of norm at most $Y$ lying outside $T$. There are at most $m$ ways of 
getting
a given element $\xi$ with (\ref{22}) by multiplying an element $\beta '$
coprime with ${\frak a}$ with an element from $\mathcal{R}$.
In other words,
there are at most $m$ principal ideals of norm $\ll X$ composed of prime
ideals of norm at most $Y$ outside $T$
which give rise
to the same $\xi$ with (\ref{22}).
Together with our lower bound $\psi_{K,T}(X,Y)/h$ for the number of
principal ideals
this implies that the number of $\xi$ with (\ref{22}) is
at least $(hm)^{-1}\psi_{K,T}(X,Y)$.
This proves Lemma 7.$\qed$
\\[0.5cm]
For functions $f(y), g(y)$ we say that $f(y)=o(g(y))$ as $y\to\infty$ 
effectively
in terms of parameters $z_1,\ldots ,z_t$ if for every $\delta >0$
there is an effectively computable constant $y_0$ depending on $\delta$,
$z_1,\ldots ,z_t$ such that $|f(y)|\leq \delta |g(y)|$ for every $y\geq y_0$. 
Then we have: 
\\[0.3cm]
\phrase{Lemma 8. } {Let $0<\alpha <1$. 
Further, let $K$ be a number field and $T$ a finite set
of prime ideals of $O_K$. Then for $Y\to\infty$ there is an $X$
such that
\begin{eqnarray}
\label{8A}
\log X&\leq &\frac{2}{1-\alpha}Y^{1-\alpha},\\[0.2cm]
\label{8B}
\frac{\psi_{K,T}(X,Y)}{X^{\alpha}}
&\geq&
\exp\big\{ \frac{1+o(1)}{1-\alpha}\cdot Y^{1-\alpha}(\log Y )^{-1}\big\}
\end{eqnarray}
where the $o$-symbol is effective in terms of $\alpha$, $K$, $T$.}
\\[0.3cm]
{\bf Proof.} Below all $o$-symbols are with respect to $Y\to\infty$ and 
effective in terms of $\alpha ,K,T$.  
Let $X=Y^u$ with $u\log u =Y^{1-\alpha}$.
Thus,
\[
u=(1+o(1))(1-\alpha )^{-1}Y^{1-\alpha}(\log Y)^{-1}
\]
and
\[
\log X =u\log Y =(1+o(1))(1-\alpha )^{-1}Y^{1-\alpha}.
\]
Note that for $Y$ sufficiently large, $X$ satisfies (\ref{8A}). Further,
$u\geq 3$.
Now by our choice of $u$ and by Theorem 5 we have
\[
\begin{array}{rl}
\displaystyle{\frac{\psi_{K,T}(X,Y)}{X^{\alpha}}}
&\geq Y^{u(1-\alpha )}
\exp\big\{-u\big({\rm log}(u\log u)~
-1~+o(1)\big)\big\}\\
&\quad\\
&\geq
\exp\big( (1+o(1))u\big)=
\exp\big\{ \frac{1+o(1)}{1-\alpha}\cdot Y^{1-\alpha}(\log Y )^{-1}\big\}
\end{array}
\]
which is (\ref{8B}).$\qed$

\section{Proofs of Theorems 1 and 2.}
\setcounter{equation}{0}
\noindent
{\bf Proof of Theorem 1.}
Constants implied by $\ll$ and $\gg$ are effective and depend
only on $n,a_1,\ldots ,a_n$ and the $o$-symbols
are always with respect to $s\to\infty$ and effective in terms of 
$n,a_1,\ldots ,a_n$. By ``sufficiently large" we mean that the quantity
under consideration exceeds some constant effectively computable in terms
of $n,a_1,\ldots, a_n$.  
We denote the cardinality of a set
$A$ by $|A|$.
\\[0.3cm]
Let $s$ be a positive integer and let $\varepsilon$ be a positive real 
number.
Put
\begin{equation}
\label{119}
\begin{array}{l}
t=[(1- \varepsilon/2) s\, ],\\[0.1cm]
Y=p_t,\,\,\, T=\{ p_1,\ldots ,p_t\}
\end{array}
\end{equation}
where
$p_i$ denotes the $i$-th prime.
Note that, by an effective version of the Prime Number Theorem,
\begin{equation}
\label{120}
Y=(1+o(1)) t \log t.
\end{equation}
We choose $X$ according to Lemma 8 with $\alpha =1/n$, $K=\Q$, $T=\emptyset$.
\\[0.2cm]
Let $\varepsilon_i = \frac{a_i}{|a_i|}$ for $i =1, \ldots , n$.
The number of $n$-tuples
$(x_1,\ldots ,x_n)$ with each $\varepsilon_i x_i$ a positive
integer of size at
most
$X$ and composed of primes at most $Y$ equals $\psi(X,Y)^n$.
Since the sum $a_1x_1+\cdots +a_nx_n$ is $\ll X$ and is a positive
rational number
with denominator $\ll 1$, there exists a
positive rational $a_0\ll X$ with denominator $\ll 1$
such that the set of tuples $(x_1,\ldots ,x_n)\in\Z^n$ with
\begin{equation}
\left.
\begin{array}{l}
a_1x_1+\cdots +a_nx_n = a_0\\[0.1cm]
1\leq \varepsilon_ix_i\leq X,\,\,
\mbox{$x_i$ is composed of primes $\leq Y$ for
$i=1,\ldots ,n$,}\\
\end{array}
\right\}
\label{4}
\end{equation}
has cardinality $\gg\psi (X,Y)^n/X$.
Let $R$ be the set of primes $p$ dividing the numerator or denominator
of $a_0$.
By the (effective) Prime Number Theorem, $|R|$ is at most
\[
(1+ o (1)) \log X/ \log_2 X.
\]
From (\ref{8A}) with $\alpha =1/n$, (\ref{120}), (\ref{119})
we infer that $|R|=o(s)$ and then from
(\ref{119}) that $|R\cup T|<s$ provided
$s$ is sufficiently large.
Let $S$ be a set of primes of cardinality $s$ containing $R\cup T$.
\\[0.3cm]
Clearly the numbers $\frac{x_i}{a_0}$ for $i = 1, \ldots , n$ are $S$-units.
Further, since $a_{i} (\frac{x_i}{a_0})$ is positive
for $i = 1,\ldots , n$, the
subsums of $a_1 x_1 + \cdots + a_n x_n$ are all non-zero.
Thus equation (\ref{1})
has $\gg \psi (X,Y)^n /X$ non-degenerate solutions in $S$-units.
By (\ref{8B}) with $\alpha =1/n$ and (\ref{120}) 
we have for $Y$ sufficiently large
\[
\begin{array}{rl}
\psi (X,Y)^n/X &\geq
\exp \Big( (1+o(1)) \frac{n^2}{n-1} Y^{1-(1/n)}
(\log Y)^{-1}\Big)\\[0.2cm]
&\geq
\exp \Big( (1+o(1)) \frac{n^2}{n-1} t^{1-(1/n)}
(\log t)^{-1/n}\Big).
\end{array}
\]
Using (\ref{119}) it follows at once that for $s$ sufficiently large,
equation (\ref{1}) has more than
$
\exp \Big( (1-\varepsilon )\frac{n^2}{n-1} s^{1-(1/n)}
(\log s)^{-1/n}\Big)$ non-degenerate
solutions in $S$-units. This proves Theorem 1.$\qed$
\\[0.5cm]
\noindent
Before proving Theorem 2 we observe that $g({\bf a},S)$ is the
smallest integer $g$ for which there exists a non-zero polynomial
$P^*\in\C [X_1,\ldots ,X_{n-1}]$ of total degree $g$ with
\begin{equation}
P^*(x_1,\ldots ,x_{n-1})=0\quad
\mbox{for every solution $(x_1,\ldots ,x_n)$ of (\ref{1}).}
\label{14}
\end{equation}
Indeed, let $P\in\C [X_1,\ldots ,X_n]$ be a polynomial of total degree
$g({\bf a},S)$ with (\ref{2}) which is not divisible by
$a_1X_1+\cdots +a_nX_n-1$. Substituting
$X_n=a_n^{-1}(1-a_1X_1-\cdots -a_{n-1}X_{n-1})$
in $P$ we get a polynomial
$P^*$ which
satisfies (\ref{14}), has total degree at most $g({\bf a},S)$,
and is not identically zero.
On the other hand, any non-zero polynomial $P^*$ with (\ref{14})
must have total degree at least $g({\bf a},S)$ since it is not
divisible by $a_1X_1+\cdots +a_nX_n-1$.
\\[0.5cm]
\noindent
{\bf Proof of Theorem 2.} Let $\varepsilon >0$.
By Theorem 1 with $n=2$
we know that there is an effectively computable positive number 
$t_1$,
which depends only on $\varepsilon$, such that
for every integer $t\geq t_1$
there is a set of primes $T$ of cardinality $t$
for which the equation
$x+y=1$ in $T$-units $x,y$ has at least
\begin{equation}
\label{111}
A(t):=\exp\big\{ (\mbox{$4-\frac{1}{2}\varepsilon$})
t^{1/2}(\log t)^{-1/2}\big\}
\end{equation}
solutions.
Fix such $t$ and $T$.
We first show by induction
that for every $n\geq 2$ the $n$-tuple
${\bf 1}_n = (1,\ldots , 1)$ satisfies $g({\bf 1}_n, T) \geq A(t)$.
\\[0.3cm]
We are done for $n=2$. Suppose $n\geq 3$, and that
our assertion holds with $n-1$ in place of $n$.
Thus $g({\bf 1}_{n-1}, T) \geq A(t)$.
Let $U$ be the set of tuples
\begin{equation}
(x_1,\ldots ,x_n)=(y_1,\ldots ,y_{n-2},y_{n-1}z_1,y_{n-1}z_2)
\label{8}
\end{equation}
where $(y_1,\ldots ,y_{n-1})$ runs through the solutions of
\begin{equation}
y_1+\cdots + y_{n-1}=1
\quad\mbox{in $T$-units $y_1,\ldots ,y_{n-1}$}
\label{10}
\end{equation}
and where $(z_1,z_2)$ runs through
the solutions of
\begin{equation}
z_1+z_2=1\quad\mbox{in $T$-units $z_1,z_2$.}
\label{11}
\end{equation}
Then from
\[
y_1 + \cdots + y_{n-2} + y_{n-1} (z_1 + z_2) = 1
\]
it follows that the tuples in $U$ satisfy
\begin{equation}
x_1+\cdots +x_n=1.
\label{9}
\end{equation}
Let $P\in \C [X_1,\ldots ,X_{n-1}]$ be a
non-zero
polynomial of total degree
$g({\bf 1}_n,T)$ such
that $P(x_1,\ldots ,x_{n-1})=0$ for every solution $(x_1,\ldots ,x_n)$ in
$T$-units of (\ref{9}).
Since the tuples in $U$ consist of $T$-units,
we have
\begin{equation}
P(y_1,\ldots ,y_{n-2},y_{n-1}z_1)=0
\label{12}
\end{equation}
for every solution $(y_1,\ldots ,y_{n-1})$ of (\ref{10}) and every solution
$(z_1,z_2)$ of (\ref{11}). Define the polynomial in $n-1$ variables
\begin{equation}
P^*(Y_1,\ldots ,Y_{n-2},Z_1)=
P(Y_1,\ldots ,Y_{n-2},\,
Z_1\cdot (1 - Y_1-\cdots - Y_{n-2})).
\label{13}
\end{equation}
Then $P^*$ is not identically zero since $P$ is not identically zero
and since the change of variables
\[
(X_1,\ldots ,X_{n-1})\mapsto (Y_1,\ldots ,Y_{n-2},\,
Z_1\cdot (1-Y_1 - \cdots - Y_{n-2}))
\]
is invertible. Now from (\ref{12}), (\ref{10}) it follows that
\begin{equation}
P^*(y_1,\ldots ,y_{n-2},z_1)=0
\label{15}
\end{equation}
for every solution $(y_1,\ldots ,y_{n-1})$ of (\ref{10}) and every solution
$(z_1,z_2)$ of (\ref{11}). We distinguish two cases.
\\[0.3cm]
{\bf Case 1.} There is a solution $(z_1,z_2)$ of (\ref{11}) such that the
polynomial\\$P^*_{z_1}(Y_1,\ldots ,Y_{n-2}):=P^*(Y_1,\ldots ,Y_{n-2},z_1)$
is not identically zero.
\\[0.2cm]
Then by (\ref{15}), $P^*_{z_1}$ is a non-zero polynomial with
$P^*_{z_1}(y_1,\ldots ,y_{n-2})=0$ for every solution
$(y_1,\ldots ,y_{n-1})$ of (\ref{10}). Hence
$P^*_{z_1}$ has total degree $\geq g({\bf 1}_{n-1},T)\geq A(t)$.
Now by (\ref{13}) this implies that the total degree
$g({\bf 1}_n,T)$ of $P$ is at least $A(t)$.
\\[0.3cm]
{\bf Case 2.} For every solution $(z_1,z_2)$ of (\ref{11}), the
polynomial $P^*_{z_1}(Y_1,\ldots ,Y_{n-2})=P^*(Y_1,\ldots ,Y_{n-2},z_1)$
is identically zero.
\\[0.2cm]
Then since (\ref{11}) has at least $A(t)$ solutions, the polynomial
$P^*$ must have degree at least $A(t)$ in the variable $Z_1$.
By (\ref{13})
this implies that $P$ has degree at least $A(t)$ in the variable $X_{n-1}$.
So again we conclude that the total degree $g({\bf 1}_n,T)$ of $P$ 
is at least $A(t)$.
This completes our induction step.
\\[0.3cm]
Now
let ${\bf a}=(a_1, \ldots , a_n)$ be an arbitrary tuple of
non-zero rational numbers and let $R$ be the set of primes dividing
the product of the numerators and denominators of $a_1,\ldots ,a_n$.
Then $|R|\ll 1$.
\\
Let $s_1$ be a positive number such that if $s$ is an integer with
$s\geq s_1$
then for
\begin{equation}
\label{112}
t :=\Big[\big(\mbox{$\frac{4-\varepsilon}{4-\varepsilon /2}$}\big)^2
\cdot s\, \Big] +1
\end{equation}
we have
\[
t\geq t_1,\quad t+|R| <s.
\]
Clearly, $s_1$ is effectively computable in terms of
$n, a_1,\ldots ,a_n, \varepsilon$. Choose $s\geq s_1$ and
let $T$ be a set of $t$ primes
with $g({\bf 1}_n,T)\geq A(t)$.
Choose any set of primes $S$ of cardinality $s$ containing
$T\cup R$.
Then since $a_1,\ldots ,a_n$ are $S$-units and by (\ref{111}), (\ref{112})
we have
\[
g({\bf a},S) =g({\bf 1}_n,S)\geq g({\bf 1}_n,T)\geq A(t)
\geq \exp\Big( (4-\varepsilon )s^{1/2}(\log s)^{-1/2}\Big).
\]
Theorem 2 follows.$\qed$

\section{Proofs of Theorems 3 and 4.}
\setcounter{equation}{0}
\noindent
We keep the notation from the previous sections.
In particular, $K$ is a number field of degree $n\geq 2$
and
$\alpha_1,\ldots ,\alpha_m$ are $\Q$-linearly independent elements of $O_K$,
where $1\leq m\leq n-1$.
Constants implied by $\ll$, $\gg$ are effectively computable in terms of
$K$, $\alpha_1,\ldots ,\alpha_m$ and the $o$-symbols
will be with respect to $s\to\infty$ and effective in terms of $K$,
$\alpha_1,\ldots ,\alpha_n$. By ``sufficiently large" we mean that the
quantity under consideration exceeds some constant effectively computable
in terms of $K,\alpha_1,\ldots ,\alpha_n$.
\\[0.2cm]
We order the rational primes $p$ by the size of
the smallest norm $p^{k_p}$ of a prime ideal dividing $(p)$.
Let $p_1, \ldots , p_s$ be the first $s$ primes in this ordering and put
$Y=p_s^{k_{p_s}}$.
By the effective version of the Chebotarev density theorem from
\cite{LaOd} (Theorems 1.3, 1.4) we have
\begin{equation}
Y=(1+o(1))c_K s \log s\, .
\label{40}
\end{equation}
\\[0.1cm]
We have to make some further preparations.
Choose $\gamma\in O_K$ with $\Q (\gamma )=K$;
then the conjugates $\gamma^{(1)},\ldots ,\gamma^{(n)}$
are distinct. Further, choose $\delta\in O_K$
which is $\Q$-linearly independent of
$\alpha_1,\ldots ,\alpha_m$. Then
there are indices
$i_0,i_1,\ldots ,i_m\in\{ 1,\ldots ,n\}$
such that
\[
\Delta :=
\left|
\begin{array}{cccc}
\alpha_1^{(i_0)}&\ldots&\alpha_m^{(i_0)}&\delta^{(i_0)}\\
\vdots&&\vdots&\vdots\\
\alpha_1^{(i_m)}&\ldots&\alpha_m^{(i_m)}&\delta^{(i_m)}
\end{array}\right|
\, \not=0.
\]
Choose a rational prime number $p$ such that $p$ is coprime with $\gamma$
and with the
differences $\gamma^{(i)}-\gamma^{(j)}$
$(1\leq i<j\leq n)$.
Further, choose another rational prime number $q$ such that $q$ is coprime
with $\delta$
and with $\Delta$.
Then by the Chinese Remainder Theorem, there is a $\beta\in O_K$
such that $\beta\equiv \gamma\pmod{p}$, $\beta\equiv\delta\pmod{q}$
and $\beta$ is coprime with $pq$. 
It is clear that $p,q,\beta$ can be determined effectively.
\\[0.3cm]
\phrase{Lemma 9. }
{For every $\xi\in O_K$ with $\xi\equiv \beta\pmod{pq}$
we have that $\Q (\xi )=K$ and that $\xi$ is $\Q$-linearly independent of
$\alpha_1,\ldots ,\alpha_m$.}
\\[0.3cm]
{\bf Proof.}
Take $\xi\in O_K$ with $\xi\equiv \beta\pmod{pq}$. Then
$\xi^{(i)}\equiv \beta^{(i)}\equiv\gamma^{(i)}\pmod{p}$ for $i=1,\ldots ,n$, 
so
\[
\xi^{(i)}-\xi^{(j)}\equiv \gamma^{(i)}-\gamma^{(j)}\not\equiv 0\pmod{p}
\]
for $1\leq i<j\leq n$, which implies that the conjugates of $\xi$
are distinct. Hence $\Q (\xi )=K$. Likewise, we have
$\xi^{(i)}\equiv \beta^{(i)}\equiv\delta^{(i)}\pmod{q}$
for $i=1,\ldots ,n$, so
\[
\left|
\begin{array}{cccc}
\alpha_1^{(i_0)}&\ldots&\alpha_m^{(i_0)}&\xi^{(i_0)}\\
\vdots&&\vdots&\vdots\\
\alpha_1^{(i_m)}&\ldots&\alpha_m^{(i_m)}&\xi^{(i_m)}
\end{array}\right|
\,\,
\equiv\Delta\not\equiv 0\pmod{q}.
\]
Hence the determinant on the left-hand side is $\not= 0$,
and therefore, $\xi$ is $\Q$-linearly independent of
$\alpha_1,\ldots ,\alpha_m$. This proves Lemma 9.$\qed$
\\[0.5cm]
{\bf Proof of Theorem 3.}
Let $V$ be the $\Q$-vector space generated by $\alpha_1,\ldots ,\alpha_m$.
Choose an integral basis $\{\omega_1,\ldots ,\omega_n\}$ of $O_K$
such that $\omega_1,\ldots ,\omega_m$ span $V$;
this can be done effectively.
Thus, every $\xi\in O_K$
can be expressed uniquely as $\xi =\sum_{j=1}^n x_j\omega_j$
with
$x_j\in\Z$. By applying Cramer's rule to
$\xi^{(i)} =\sum_{j=1}^n x_j\omega_j^{(i)}$ $(i=1,\ldots ,n)$
and using the fact that det$\,(\omega_j^{(i)})\not=0$ we get
\[
\max_{j=1,\ldots ,n} |x_j| \,\ll \, \max_{i=1,\ldots ,n} |\xi^{(i)}|.
\]
We combine this with Lemma 7. Choose $X>Y$.
Since by our construction, $\beta$ is coprime
with $pq$, it follows that the set of $\xi\in O_K$
with
\[
\begin{array}{l}
\xi =\sum_{j=1}^n x_j\omega_j,
\quad\mbox{$x_j\in\Z,\quad |x_j|\ll X^{1/n}$ for $j=1,\ldots ,n$,}
\\[0.1cm]
\xi\equiv\beta\pmod{pq},
\\[0.1cm]
\mbox{$(\xi )$ composed of prime ideals of norm $\leq Y$}
\end{array}
\]
has cardinality $\gg \psi_{K,T}(X,Y)$,
where $T$ is the set of prime ideals dividing $(pq)$. Consequently, there is 
a number
\[
\kappa = \sum_{j=m+1}^n y_j\omega_j\quad
\mbox{with $y_j\in\Z$, $|y_j|\ll X^{1/n}$ for $j=m+1,\ldots ,n$}
\]
such that the set of $\xi\in O_K$ with
\begin{equation}
\left.
\begin{array}{l}
\xi =\kappa +\sum_{j=1}^m x_j\omega_j,
\quad\mbox{$x_j\in\Z,\,\, |x_j|\ll X^{1/n}$ for $j=1,\ldots ,m$,}
\\[0.1cm]
\xi\equiv\beta\pmod{pq},
\\[0.1cm]
\mbox{$(\xi )$ composed of prime ideals of norm $\leq Y$}
\end{array}
\right\}
\label{43}
\end{equation}
has cardinality $\gg \psi_{K,T}(X,Y)/X^{1-(m/n)}$.
\\[0.3cm]
Pick $\xi_0$ satisfying (\ref{43}). Then by Lemma 9,
$\xi_0$ is an algebraic integer such that $\Q (\xi_0 )=K$ and
$\xi_0$ is $\Q$-linearly independent of $\alpha_1,\ldots ,\alpha_m$.
Since $\omega_1,\ldots ,\omega_m$ span the same $\Q$-vector space
as $\alpha_1,\ldots ,\alpha_m$, there is a positive rational integer $d$
such that the $\Z$-module generated by $d\omega_1,\ldots ,d\omega_m$
is contained in the $\Z$-module generated by $\alpha_1,\ldots ,\alpha_m$.
Put $\alpha_0:=d\xi_0$; then $\alpha_0$ satisfies (\ref{33}).
\\[0.3cm]
We have $\xi_0 =\kappa +\sum_{j=1}^m y_j\omega_j$ with $y_j\in\Z$,
$|y_j|\ll X^{1/n}$ for $j=1,\ldots ,m$.
If for $\xi$ as in (\ref{43}) we write $x'_j=x_j-y_j$ $(j=1,\ldots ,m)$,
we get
\[
\xi =\xi_0 +\sum_{j=1}^m x'_j\omega_j\quad
\mbox{with $x'_j\in\Z,\,\,\,|x'_j|\ll X^{1/n}$ for $j=1,\ldots ,m$}
\]
(where we have enlarged the constant implied by $\ll$). By expressing
$d\omega_1,\ldots ,d\omega_m$ as linear combinations of
$\alpha_1,\ldots ,\alpha_m$ with coefficients in $\Z$ we may express
$d\xi$ with $\xi$ satisfying (\ref{43}) as
\begin{equation}
d\xi =
\alpha_0+\sum_{j=1}^m x''_j\alpha_j\quad
\mbox{with $x''_j\in\Z,\,\,\,|x''_j|\ll X^{1/n}$ for $j=1,\ldots ,m$}
\label{44}
\end{equation}
(again after enlarging the constant implied by $\ll$).
Assuming, as we may, that $d$ is
composed of prime ideals of norm at most $Y$, we have for $\xi$ with
(\ref{43})
that $(d\xi )$ is composed of prime ideals of norm at most $Y$.
Hence $|N_{K/\Q}(d\xi )|$ is composed of $p_1,\ldots ,p_s$.
To simplify notation we write $x_j$ instead of $x''_j$.
Recalling that the set of elements
with (\ref{43})
has cardinality $\gg \psi_{K,T}(X,Y)/X^{1-(m/n)}$
and that $d\xi$ with $\xi$ as in (\ref{43}) can be expressed as (\ref{44}),
we obtain that the set of tuples $(x_1,\ldots ,x_m)\in\Z^m$
with
\begin{equation}
\left.
\begin{array}{l}
|N_{K/\Q }(\alpha_0+x_1\alpha_1+\cdots +x_m\alpha_m)|
=p_1^{z_1}\cdots p_s^{z_s}\,\,
\mbox{for certain $z_1,\ldots ,z_s\in\Z$,}\\[0.2cm]
\mbox{$|x_j|\ll X^{1/n}$ for $j=1,\ldots ,m$}
\end{array}
\right\}
\label{42}
\end{equation}
has cardinality $\gg \psi_{K,T}(X,Y)/X^{1-(m/n)}$.
\\[0.3cm]
We have already observed that $Y\to\infty$ as $s \rightarrow \infty$. 
Further,
from Lemma 8 with $\alpha =1-(m/n)$ and from (\ref{40}) it
follows that for arbitrarily large $Y$ there is an $X$ with
\[
\begin{array}{rl}
\psi_{K,T}(X,Y)/X^{1-(m/n)}
&\geq \exp \big\{ (1+o(1))\frac{n}{m}\cdot Y^{m/n}(\log Y)^{-1}\big\}\\
&\quad\\
&\geq \exp \big\{ (1+o(1))\frac{n}{m}\cdot (c_Ks)^{m/n}(\log
s)^{(m/n)-1}\big\}.
\end{array}
\]
Theorem 3 now follows directly.$\qed$
\\[0.5cm]
{\bf Proof of Theorem 4.} In the proof of Theorem 3 we have shown
that for every sufficiently large $Y$ and every $X>Y$
there is an $\alpha_0$ with (\ref{33}), such that the set of
tuples ${\bf x}=(x_1,\ldots ,x_m)\in\Z^m$ with (\ref{42})
has cardinality
$\gg \psi_{K,T}(X,Y)/X^{1-(m/n)}.$
\\[0.3cm]
Let $S=\{ p_1,\ldots ,p_s\}$.
Let $P\in\C [X_1,\ldots ,X_m]$
be a non-trivial polynomial of total degree $g=g(\ba ,S)$
such that for each solution $(x_1,\ldots ,x_m,\, z_1,\ldots ,z_s)$
of (\ref{31}) we have $P(x_1,\ldots ,x_m)=0$. This implies in particular
that $P({\bf x})=0$ for each tuple ${\bf x}$ with (\ref{42}).
Now since the tuples with (\ref{42}) have $|x_j|\ll X^{1/n}$
for $j=1,\ldots ,m$, we have by Lemma 6 that there are
$\ll g\cdot (X^{1/n})^{m-1}$ vectors $(x_1,\ldots ,x_m)$ with (\ref{42}).
Together with our lower bound $\gg \psi_{K,T}(X,Y)/X^{1-(m/n)}$
for the number of these vectors, this gives
\[
g\cdot X^{(m-1)/n}\gg \psi_{K,T}(X,Y)/X^{1-(m/n)}
\]
or equivalently
\[
g\gg \psi_{K,T}(X,Y)/X^{1-(1/n)}.
\]
Again $Y$ goes to infinity with $s$. Further,
by Lemma 8 with $\alpha =1-(1/n)$ and (\ref{40})
we have that for $Y\to\infty$
there is an $X$ with
\[
\begin{array}{rl}
\psi_{K,T}(X,Y)/X^{1-(1/n)}
&\geq \exp \big\{ (1+o(1))n\cdot Y^{1/n}(\log Y)^{-1}\big\}\\
&\quad\\
&\geq \exp \big\{ (1+o(1))n\cdot (c_Ks)^{1/n}(\log s)^{(1/n)-1}\big\}.
\end{array}
\]
This proves Theorem 4.$\qed$

\newpage
{\bf Adresses of the authors:}
\\[0.3cm]
J.-H. Evertse and R. Tijdeman:\\
Universiteit Leiden, Mathematisch Instituut,\\
Postbus 9512, 2300 RA Leiden, The Netherlands\\
email evertse@math.leidenuniv.nl, tijdeman@math.leidenuniv.nl
\\[0.3cm]
P. Moree:\\
Universiteit van Amsterdam, KdV Instituut,\\
Plantage Muidergracht 24, 1018 TV Amsterdam, The Netherlands\\
email moree@science.uva.nl
\\[0.3cm]
C.L. Stewart:\\
University of Waterloo, Department of Pure Mathematics,\\
Waterloo, Ontario, Canada, N2L 3G1\\
email cstewart@watserv1.uwaterloo.ca

\end{document}